\documentclass{article}

     \PassOptionsToPackage{numbers, compress}{natbib}

\usepackage[final]{neurips_2023_ml4ps}

\newcommand{\paren}[1]{\ensuremath{\left( #1 \right)}}

\newcommand{\curlyb}[1]{\ensuremath{\left\{ #1 \right\} }}

\newcommand\thetabold{{\ensuremath{\boldsymbol{\theta}}}}

\usepackage{bm} 
\usepackage{amsfonts} 

\newcommand{\norm}[1]{\ensuremath{\left\| #1 \right\|}}

\newcommand{\pder}[2]{\ensuremath{\frac{\partial #1}{\partial #2}}}

\newcommand{\Dcal}{\ensuremath{\mathcal{D}}}

\newcommand{\Lcal}{\ensuremath{\mathcal{L}}}

\newcommand{\Ocal}{\ensuremath{\mathcal{O}}}

\newcommand{\Rbb}{\ensuremath{\mathbb{R} }}

\newcommand\Dbm{{\ensuremath{\bm{D}}}}
\newcommand\Ebm{{\ensuremath{\bm{E}}}}

\newcommand\Ibm{{\ensuremath{\bm{I}}}}

\newcommand\Ubm{{\ensuremath{\bm{U}}}}
\newcommand\Vbm{{\ensuremath{\bm{V}}}}
\newcommand\Wbm{{\ensuremath{\bm{W}}}}

\newcommand\bbm{{\ensuremath{\bm{b}}}}

\newcommand\fbm{{\ensuremath{\bm{f}}}}

\newcommand\hbm{{\ensuremath{\bm{h}}}}

\newcommand\rbm{{\ensuremath{\bm{r}}}}
\newcommand\sbm{{\ensuremath{\bm{s}}}}

\newcommand\ubm{{\ensuremath{\bm{u}}}}
\newcommand\vbm{{\ensuremath{\bm{v}}}}
\newcommand\wbm{{\ensuremath{\bm{w}}}}
\newcommand\xbm{{\ensuremath{\bm{x}}}}

\newcommand\zbm{{\ensuremath{\bm{z}}}}

\newcommand\mubold{{\ensuremath{\boldsymbol{\mu}}}}

\newcommand\phibold{{\ensuremath{\boldsymbol{\phi}}}}

\newcommand\psibold{{\ensuremath{\boldsymbol{\psi}}}}

\newcommand\sigmabold{{\ensuremath{\boldsymbol{\sigma}}}}

\newcommand\zerobold{\ensuremath{\mathbf{0}}}

\usepackage[utf8]{inputenc} 
\usepackage[T1]{fontenc}    
\usepackage{hyperref}       
\usepackage{url}            
\usepackage{booktabs}       
\usepackage{amsfonts}       
\usepackage{nicefrac}       
\usepackage{microtype}      
\usepackage{xcolor}         
\usepackage{color}
\usepackage{svg}
\definecolor{green2}{rgb}{0, 0.5, 0}

\usepackage{float}
\usepackage{lscape}

\title{Reduced-order modeling for parameterized PDEs\\ via implicit neural representations}

\author{%
  Tianshu Wen\\
  Aerospace and Mechanical Engineering\\
  University of Notre Dame\\
  \texttt{twen2@nd.edu} \\
  \And
  Kookjin Lee \\
  Computing and Augmented Intelligence \\
  Arizona State University \\
  \texttt{kookjin.lee@asu.edu} \\
  \AND
  Youngsoo Choi \\
  Center for Applied Scientific Computing \\
  Lawrence Livermore National Laboratory \\
  \texttt{choi15@llnl.gov} \\
}

\usepackage{bm}
\usepackage{physics}
\usepackage{mathtools}
\mathtoolsset{showonlyrefs}
\usepackage{amssymb,amsthm}
\usepackage{algorithmic}
\usepackage[linesnumbered,ruled,noend]{algorithm2e} 

\usepackage{multirow}
\usepackage{wrapfig}

\usepackage{caption}
\usepackage{subcaption}

\usepackage{tabularx,ragged2e,booktabs,lipsum}

\newcommand{\Transpose}{^{\mathsf{T}}}

\begin{document}

\maketitle

\begin{abstract}
We present a new data-driven reduced-order modeling approach to efficiently solve parametrized partial differential equations (PDEs) for many-query problems. 
This work is inspired by the concept of implicit neural representation (INR), which models physics signals in a continuous manner and independent of spatial/temporal discretization.
The proposed framework encodes PDE and utilizes a parametrized neural ODE (PNODE) to learn latent dynamics characterized by multiple PDE parameters. PNODE can be inferred by a hypernetwork to reduce the potential difficulties in learning PNODE due to a complex multilayer perceptron (MLP).
The framework uses an INR to decode the latent dynamics and reconstruct accurate PDE solutions. Further, a physics-informed loss is also introduced to correct the prediction of unseen parameter instances. Incorporating the physics-informed loss also enables the model to be fine-tuned in an unsupervised manner on unseen PDE parameters.
A numerical experiment is performed on a two-dimensional Burgers equation with a large variation of PDE parameters. We evaluate the proposed method at a large Reynolds number and obtain up to speedup of $\Ocal(10^3)$ and $ \sim 1\%$ relative error to the ground truth values. 
\end{abstract}
\section{Introduction}\label{sec:intro}
Numerical simulations are broadly used in almost every branch of science and engineering, such as climate modeling, product design, risk prediction, etc. However, high-fidelity simulations remain challenging in practice due to the high dimensionality and complicated physics patterns, especially for many-query problems. To circumvent the computational cost in PDE simulations, a slew of surrogate modeling techniques have been developed in many forms. 

One particular approach is the projection-based reduced-order model (ROM) which has succeeded in many applications \cite{zahr_progressive_2015, arian_trust-region_2000, qian_certified_2017, choi2020gradient, yano_globally_2021, mcbane2022stress, wen_globally_2023, carlberg2018conservative, cheung2023datascarce, lauzon2022s, tsai2023accelerating, kim2021efficient, copeland2022reduced, huhn2023parametric, choi2020sns, cheung2023local, mcbane2021component, hoang2021domain, choi2019space}. Many traditional ROMs utilize proper orthogonal decomposition (POD) based linear projection to reconstruct accurate approximations. However, due to the limitation of the linear projections, the Kolmogorov $n$-width problems remain an obstacle to those linear ROMs. 

In recent years, neural network based ROM approaches have been developed \cite{swischuk_projection-based_2019, lee_model_2020, lee2021deep, kim2020efficient, kim2022fast, fries_lasdi_2022, kadeethum_non-intrusive_2022, he_glasdi_2023, barnett_neural-network-augmented_2023}. The capability of machine learning (ML) to make universal approximations allows such ROMs to approximate PDE states in a nonlinear manifold.
ROM with neural networks utilizes encoder-decoder-type structures in which they learn the latent dynamics in low-dimensional latent space and then decode the reduced states to reconstruct the PDE states. Among various ML-based approaches, implicit neural representation (INR) type decoders are gaining traction  \cite{chen2022crom, wan2023evolve, yin_continuous_2023} due to their flexibility; INRs decode the latent states at the continuous level, meaning that it has the capability to extrapolate arbitrary spatial/temporal positions and is independent of the resolution of the discretization. 

In this work, we extend the state-of-the-art INR-based ROM, Dynamics-aware Implicit Neural Representation (DINo) \cite{yin_continuous_2023}, to learn the surrogate models for the solutions of parameterized PDEs. DINo compresses/decompresses high-dimensional PDE states into low-dimensional embeddings via INR-based auto-decoding/decoding and describes the temporal evolution of embeddings via neural ordinary differential equations (NODEs) \cite{chen2018neural}. However, NODEs are known to be restrictive in modeling parameterized dynamics and, thus, we extend the latent dynamics of DINo to the parameterized NODEs (PNODEs) in \cite{lee_parameterized_2021} and its hypernetwork-based alternative (HyperPNODE).
Furthermore, we introduce a fine-tuning via minimizing a physics-informed (PI) loss \cite{raissi2019physics} to correct prediction for unseen data. This "unsupervised" loss helps the decoder in producing better predictions on out-of-distribution PDE parameters. The main contributions of this work are summarized as:
\begin{itemize}
    \item a data-driven reduced-order modeling framework that is based on INR to generate spatially/temporally continuous predictions,
    \item extending the non-parametric approach (DINO) to parametrized PDE with PNODE and HyperPNODE,
    \item inference corrected by the physics-informed loss,
    \item demonstrated performance on parametrized PDE and unseen parameter instances.
    \item fine-tune the pre-trained model to improve the accuracy on unseen parameter instances.
\end{itemize}

\section{Model order reduction}\label{sec:rom}
\paragraph{Full-order model}
We consider a parameterized system of semi-discrete PDE or a system of ODEs, which we will refer to as a full order model (FOM):
\begin{equation}
\label{eqn:fom}
    \pder{\ubm}{t} = \fbm(\ubm, t;\mubold), \quad \ubm(0;\mubold) = \ubm_0(\mubold),
\end{equation}
where $\fbm: \Rbb^{N_\ubm} \times [0,T] \times \Dcal \rightarrow \Rbb^{N_\ubm}$ denotes the velocity, or the spatially discretized PDE, $\ubm(t;\mubold)$, $\ubm: [0,T] \times \Dcal \rightarrow \Rbb^{N_\ubm}$ is the PDE state and implicitly defined as the solution to the system of ODEs, $\mubold \in \Dcal \subset \Rbb^{N_\mubold}$ denotes a collection of PDE parameters, and $t$ is the time from 0 to the final time $T \in \Rbb_+$. Finally, the initial condition is specified by $\ubm_0(\mubold)$, $\ubm_0: \Dcal \rightarrow \Rbb^{N_\ubm}$ in the parametrized setting.
Solving Eq.~\eqref{eqn:fom} can be computationally expensive due to high degrees of freedom (typically $\sim 10^7$) in practical problems in computational physics. 

\paragraph{Reduced-order model}
ROMs migrate the computational cost of FOMs via latent-state evolutions in a low-dimensional manifold and make predictions on the FOM solutions via a decoder, 
a mapping from the latent states to the full-order states, $\Dbm: \Rbb^k \mapsto \Rbb^{N_\ubm}$, which follows the framework of the latent space dynamics identification (LaSDI) \cite{fries_lasdi_2022, he_glasdi_2023, bonneville2024gplasdi, tran2023weakform}. For parameterizing the latent dynamics,  we largely follow the approaches considered in \cite{yin_continuous_2023,wan2023evolve}, namely, employing NODEs \cite{chen2018neural}:
\begin{equation}\label{eqn:node}
    \pder{\hat\ubm}{t} = \hat \fbm_{\psibold} (\hat\ubm, t; \mubold, \psibold), \quad \hat\ubm(0; \mubold) = \hat\ubm_0(\mubold),
\end{equation}
where $\hat \fbm_\psibold(\cdot, \cdot; \cdot, \psibold): \Rbb^k \times [0,T] \times \Dcal \rightarrow \Rbb^k$ is a velocity function that defines the dynamics of the latent states over time, acting as a surrogate model for FOMs with $k \ll N_{\ubm}$, and $\psibold$ consists of neural network weights. Also, $\hat\ubm(t;\mubold)$, $\hat\ubm: [0,T] \times \Dcal \rightarrow \Rbb^k$ are the reduced states, and $\hat\ubm_0(\mubold)$, $\hat\ubm_0: \Dcal \rightarrow \Rbb^k$ denotes the reduced initial condition. In the ROM setting, nonlinear mapping and latent-dynamics models are trained to generate accurate approximations to the full-order model solution, i.e., $\Dbm(\hat\ubm) \approx \ubm$. While effective in learning certain classes of PDE solutions \cite{yin_continuous_2023,wan2023evolve}, these frameworks based on NODEs naturally fail to build surrogate models for dynamical systems, where the input PDE parameters (e.g., Reynolds number) can change the model dynamics. This is because of NODEs' limited expressivity (Eq.~\eqref{eqn:node}), which fails to capture the model dynamics based on the input PDE parameters. For notational simplicity, we use $\hat\ubm$ to denote reduced states or latent states for all reduced settings including NODE, PNODE, and HyperPNODE. 


\paragraph{Learning latent-dynamics using PNODE}

To address such limitation, a parameterized alternative of NODEs, named PNODEs   \cite{lee_parameterized_2021}, has extended NODEs by introducing PDE parameters $\mubold$ into the dynamics:
\begin{equation}\label{eqn:pnode}
    \pder{\hat\ubm}{t} = \hat \fbm_{\phibold} (\hat\ubm, t; \mubold, \phibold), \quad \hat\ubm(0; \mubold) = \hat\ubm_0(\mubold)
\end{equation}
where $\hat \fbm_\phibold(\cdot, \cdot; \cdot, \phibold): \Rbb^k \times [0,T] \times \Dcal \rightarrow \Rbb^k$ is the parametrized velocity function. PNODEs have been demonstrated to present multiple trajectories characterized by each PDE parameter instance. Since $\fbm_\phibold$ is an MLP, $\phibold = \{\Wbm_l,\bbm_l \}_{l=1}^{L}$ includes the network's weights and biases for $L$ layers. Inferring high-dimensional $\Wbm_l$ can be challenging in complex MLPs. 

\paragraph{Learning latent-dynamics using HyperPNODE}
To overcome the potential high complexity in MLP, we propose a low-rank factored representation of model weights and use hypernetwork to infer some parts of the model parameters. In this low-rank setting, the network module i) takes the initial latent state $\hat\ubm_0$, a time interval $t$, and PDE parameter $\mubold$ as inputs, ii) mimics the singular value decomposition (SVD), i.e., a linear combination of low-rank matrices and decomposes the weights of the MLP $\Wbm^l(\mubold) = \Sigma^r_{i=1} {s_i}^l (\mubold) \ubm^l_i \vbm^l_i {}\Transpose$, where $s_i(\mubold)$ comes from the hypernetwork, and iii) predicts a sequence of latent states $\hat\ubm(T; \mubold)$ for the selected time interval. Note that $\ubm_i$ in the SVD-like decomposition differs from the PDE states in Eq. (\ref{eqn:fom}).

\begin{figure}
    \centering
    \includegraphics[width=\textwidth]{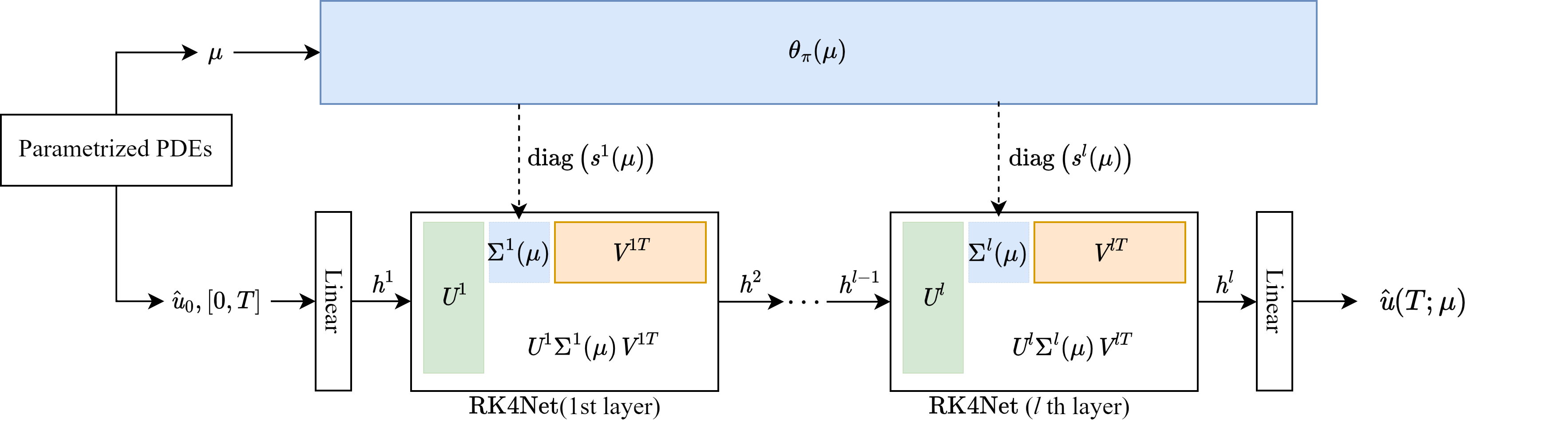}
    \caption{The architecture of HyperPNODE consisting of i) a hypernetwork to generate model parameters (i.e., diagonal elements) and ii) an RK4Net as an integrator for time interval $[0, T]$}
    \label{fig:hypernet}
\end{figure}

Next, we discuss the detailed architecture of the hypernetwork setting.
We denote this approach as a hypernetwork-based PNODE (HyperPNODE):
\begin{equation}
    \pder{\hat\ubm}{t} = \hat \fbm_{\thetabold} (\hat\ubm, t; \mubold, \thetabold), \quad \hat\ubm(0; \mubold) = \hat\ubm_0(\mubold)
\end{equation}
where $\hat \fbm_\thetabold(\cdot, \cdot; \cdot, \thetabold): \Rbb^k \times [0,T] \times \Dcal \rightarrow \Rbb^k$ is the parametrized velocity function inferred by a hypernetwork $\thetabold \coloneqq \thetabold_{\pi}(\pmb \mu)$ which only depends on the PDE parameter $\mubold$ and is learned via training. The hypernetwork $\thetabold$ infers only the diagonal entries $\{\sbm^l (\mubold)\}_{l=2}^{L-1} = \thetabold$.
Then, for a selected time interval $[0, T]$, the predicted sequence of latent states can be represented as
\begin{equation}
    \hat\ubm(T; \mubold) = \hat\ubm_0(\mubold) + \int_0^T \hat \fbm_{\thetabold} (\hat\ubm, t; \mubold, \thetabold) dt = \sigmabold(\hat \fbm_{\thetabold} (\hat\ubm, t; \mubold, \thetabold)),
\end{equation}
where $\sigmabold$ is an ODE solver. Unlike the black box integrator used in NODE and PNODE, in the proposed work, we employ the RK4Net in \cite{krishnapriyan_learning_2022} to integrate the latent dynamics. 
Finally, the internal layers of HyperPNODE can be described as:
\begin{equation}\label{eqn:low_rank_mlp}
    \begin{aligned}
        \hbm^{1} &= \sigmabold \paren{\Wbm^{0} \hbm^0 + \bbm^{0}}, \\
        \bm \hbm^l &=  \sigmabold \paren{\Ubm^{l-1} \text{diag}(\sbm^{l-1}(\mubold)) \Vbm^{(l-1)}{}\Transpose{} \hbm^{l-1} + \bbm^{l-1}}, \qquad l=2,\ldots,L-1,\\
        \hat\ubm(T; \mubold) &= \sigmabold \paren{W^L \hbm^L + \bbm^L},
    \end{aligned}
\end{equation}
where $\Ubm \in \Rbb^{n \times r}$ and $\Vbm \in \Rbb^{n \times r}$ are orthogonal matrices, and both have a size of $n$ features by $r$ user-defined rank. 
By choosing a low rank $r < k$,  $\hat \fbm_\thetabold$ is relatively efficient to learn. Figure. \ref{fig:hypernet} summarizes the overall architecture of the HyperPNODE.

\paragraph{Decoder via INR} We decode the latent states using the FourierNet~\cite{fathony_multiplicative_2021}, following \cite{yin_continuous_2023}. It outperforms the other INR architectures in cases without prior knowledge of the PDE and generalizing multiple trajectories. In this setting, the decoder $\Dbm_{\thetabold_\mathrm{dec}}$ parametrized by $\thetabold_\mathrm{dec}$ is defined as:
\begin{equation}\label{eqn:decoder}
    \tilde\ubm = \Dbm_{\thetabold_\mathrm{dec}}(\hat\ubm; \xbm, \thetabold_\mathrm{dec}),
\end{equation}
where $\tilde\ubm \in \Rbb^{N_\ubm}$ are the approximated high-dimensional states for the time interval t. The FourierNet predicts the PDE states at all nodal positions in $\Omega$ and depends on the latent states $\hat\ubm$. 

\paragraph{Training} In this work, we use an auto-decoder in \cite{park_deepsdf_2019} to encode the high-dimensional states into latent space. After the networks are defined, the forward pass of the proposed framework can be summarized as:
\begin{enumerate}
    \item encode the reduced initial state from given initial condition: $\hat \ubm_0(\mubold) = \Ebm_{\thetabold_\mathrm{enc}}(\ubm_0(\mubold);\thetabold_\mathrm{enc})$
    \item learn the latent states $\hat \ubm$ using PNODE or HyperPNODE
    \item decode the collected latent states and use INR to reconstruct high-dimensional states.
    \item compute the loss function that consists of data-matching loss and physics-informed loss (PDE residual loss, initial condition loss, and boundary condition loss) \cite{raissi2019physics}
    \begin{equation*}
        \Lcal = \alpha_1 \underbrace{\mathrm{MSE}(\ubm, \tilde\ubm)}_{\text{data-matching}} 
        + \alpha_2 \underbrace{\mathrm{MSE}(\rbm(\tilde\ubm), \zerobold)}_{\text{PDE residual}}
        + \alpha_3 \underbrace{\mathrm{MSE}(\tilde\ubm_0, \ubm_0)}_{\text{initial condition}}
        + \alpha_4 \underbrace{\mathrm{MSE}(\tilde\ubm_{\partial \Omega}, \ubm_{\partial \Omega})}_{\text{boundary condition}}
    \end{equation*}
    where $\rbm$ is the residual of the governing equation and $\alpha_i$ are user-defined scaling number. 
\end{enumerate}
For the training loss, the main loss is the data-matching loss and PDE residual loss can be used as an option. When the PDE residual loss is opted in for training, the methods are denoted by using a name with a suffix, ``+PI''. Also, when HyperPNODE is chosen, an additional orthogonality penalty is minimized: $\rho_1\| \Ubm\Transpose{}\Ubm - \Ibm \| + \rho_2\| \Vbm\Transpose{}\Vbm - \Ibm \|$, where $\rho_i$ are the penalty weights. 

\paragraph{Fine-tuning INR} As in many other ROM approaches, the proposed method shares the same limitation, performance degradation for unseen PDE parameters (i.e., out-of-distribution samples). To alleviate this limitation, we fine-tune the trained model on a target test PDE parameter in an unsupervised manner by minimizing only the PDE residual loss and selected part of latent dynamics. 
To this end, we fine-tune the INR with the loss:
\begin{equation*}\label{eqn:fine_tune_loss}
    \Lcal_{\text{ft}} = \beta_1 \mathrm{MSE}(\rbm(\tilde\ubm), \zerobold) + \beta_2 \mathrm{MSE}(\tilde\ubm_0, \ubm_0) + \beta_3 \mathrm{MSE}(\tilde\ubm_{\partial \Omega}, \ubm_{\partial \Omega}),
\end{equation*}
where $\beta_i$ are penalty weights for PDE residuals, initial condition, and boundary conditions, respectively. Additionally, we fine-tune the hypernetwork since it has a manageable number of parameters.

\section{Numerical Results and Discussion}
In this section, we apply the proposed INR-ROM method to solve a parametrized computational physics problem adapted from \cite{fries_lasdi_2022}. We demonstrate the effectiveness of the proposed parametric method by comparing the performance with using NODE.
\paragraph{2D Burgers}
The problem we consider is solving a two-dimensional, parametrized Burgers' equation:
\begin{equation}\label{eqn:burgers}
    \pder{\vbm}{t} = -\vbm \cdot \nabla \vbm + \frac{1}{\mubold} \Delta \vbm , \quad \Omega = [-3, 3] \times [-3, 3], \quad t \in [0,1],
\end{equation}
with the boundary condition $\vbm(\xbm, t;\mubold) = 0$ on $\partial \Omega$, an initial condition $\vbm(\xbm, 0;\mubold) = 0.8e^{-\frac{\norm{\xbm}^2}{1.02}}$ and $\mubold$ is the Reynolds number used as the PDE parameter. The computational domain has a uniform spatial discretization with $64\times64$ grid. The full-order model utilizes the implicit backward Euler time integrator with a uniform time step of $\mathit{\Delta} t = 1/1000$.
In this problem, the training set $\mubold_\mathrm{train} = \curlyb{30, 50, 100, 500, 1000, 2000, 5000, 10000, 30000, 50000}$, the testing set $\mubold_\mathrm{test} = \curlyb{20, 300, 20000, 60000}$ are chosen to ensure a large variance on $\mubold$, which is challenging for the physics-informed neural networks. For each model, we set the dimension of the latent states to be $50$ for each solution component (two components $\wbm$ and $\zbm$ in 2D Burgers equation), and ODE MLPs (NODE, PNODE, and RK4Net) to have 3 layers with 256 units per layer. In HyperPNODE, the hypernetwork has only 1 layer with 50 units, and the rank $r = 50$ is chosen empirically for reasonable training accuracy. Finally, the learning rate of the INR is $0.01$ and $0.001$ for all other optimizers. We also note that all results shown in this section are the first component $\wbm$.

We first investigate the performance of different approaches including NODE, PNODE, and HyperPNODE as well as their PI variants on the training set. Table. \ref{tab:train} lists the number of parameters and point-wise relative error to the ground truth values from FOM. Both hypernetworks-based approaches outperform other methods with fewer trainable parameters.
Figure. \ref{fig:train} compares the point-wise difference of the PDE states from selected models on the training point $Re = 500$. In this parametric setting, NODE fails all cases to predict the ground truth because it minimizes an MSE loss without recording parameter footprints so that it overfits an averaged $\mubold$ during the training. Hence, we only discuss the accuracy of PNODE and HyperPNODE in the rest of this discussion.
After adding physics-informed loss as shown in \ref{fig:train_pinn}, we observe 
a slightly larger error for each model. However, we observe that a noisy region is pushed towards the location near discontinuity. This indicates a potential improvement in the physics-informed variants. We cut off the training at 50000 epochs, but a model with PI loss may require further training for better results. 

Next, we test the pre-trained model on the test dataset that includes interpolated and extrapolated test values. Table. \ref{tab:test} shows similar results in the training results that adding PI loss can slightly enlarge the error. Meanwhile, the HyperPNODE models have large test errors in low Re cases.  
Figures \ref{fig:test} and \ref{fig:test_pinn} show similar behaviors as in the training results. Adding PI loss increases the maximum error, but it also pushes the error to a reasonable region. 

\begin{figure}[H]
\centering
\begin{minipage}[t]{0.32\textwidth}
\centering
\includegraphics[width=1\textwidth]{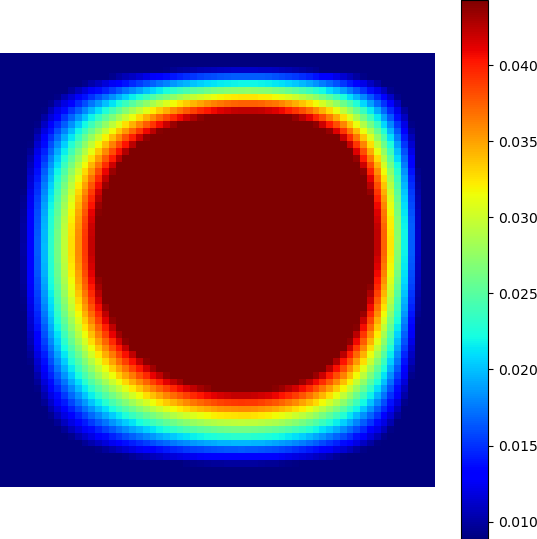}
\end{minipage}
\hfill
\begin{minipage}[t]{0.32\textwidth}
\centering
\includegraphics[width=\textwidth]{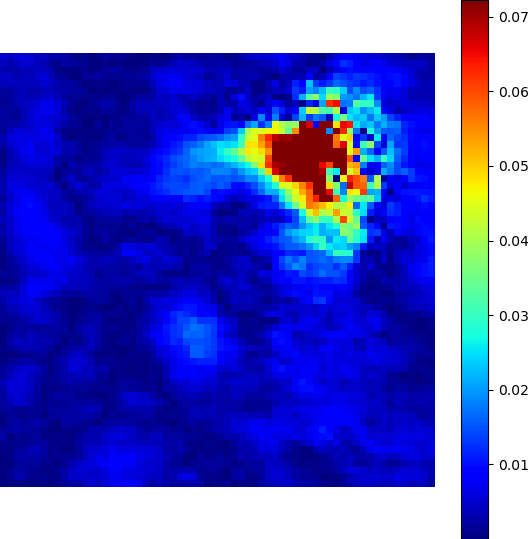}
\end{minipage}
\hfill
\begin{minipage}[t]{0.32\textwidth}
\centering
\includegraphics[width=\textwidth]{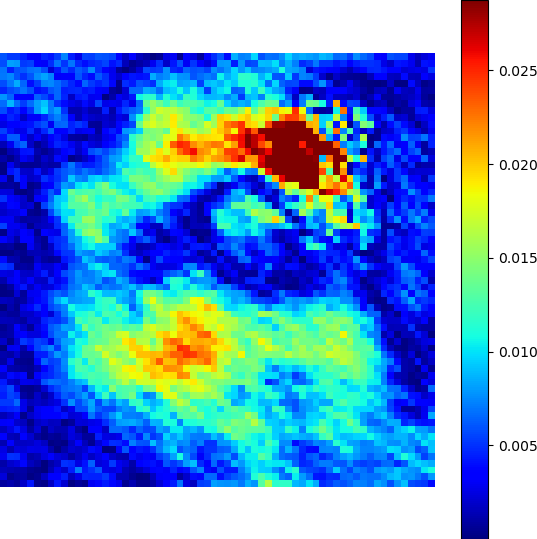}
\end{minipage}
\caption{\small Comparison of point-wise $\norm{\ubm_{\text{ts}} - \hat\ubm_{\text{ts}}}$ at $T=1$ for selected $Re = 20$ in the testing set. From left to right: reference full order model (FOM) solution, point-wise error between HyperPNODE and FOM, and point-wise error between fine-tuned HyperPNODE and FOM. 
}
\label{fig:fine_tune}
\end{figure}

Finally, we fine-tune the model at the test data point, $Re = 20$, which introduces the largest relative error. Figure. \ref{fig:fine_tune} shows the fine-tuned accuracy by HyperPNODE. The relative test error at $Re=20$ is reduced from 1.234e-01 to 1.005e-01, which is smaller than the result of HyperPNODE without PI loss. In this case, any model without incorporating physics information in the training stage fails this fine-tuning step since it does not record any physics information.
\paragraph{Discussion} 
In this work, we propose a new data-driven reduced-order model inspired by DINo for parametric configurations by using parametrized Neural ODEs (PNODE and HyperPNODE). Both PNODE and HyperPNODE successfully reconstruct the FOM solutions and provide reasonable solutions on unseen PDE parameters. With all variances of the proposed framework, we obtain $\sim 1\%$ averaged relative error on the training set, and $\sim 4\%$ averaged error on the test set. Different models with their PINN variants have a similar range of speedups from $\Ocal(10^2)$ to $\Ocal(10^3)$. We also employ an unsupervised fine-tuning method to improve the worst scenario of the HyperPNODE+PI loss model by reducing $\sim 2\%$ of the error. 

\paragraph{Future direction} There are potential improvements to the current PINN variants, including tailoring the FouriNet for spatial-temporal inputs, additional training stage, and reducing the overfitting issue in HyerPNODE models.

\section{Broader impact}
This paper introduces a new data-driven parametrized reduced-order-modeling technique using INR as a decoder to reconstruct data at a continuous level. The proposed INR-ROM framework is expected to have broad impacts on the computational science
community and application potentials in a wide range of engineering and scientific domains. There is no negative consequence on ethics and society in this work.

\begin{ack}
This work was performed under the auspices of the U.S. Department of Energy (DOE), by Lawrence Livermore National Laboratory (LLNL) under Contract No. DE-AC52–07NA27344.
Y.\ Choi was supported for this work by the U.S. Department of Energy, Office of Science, Office of Advanced Scientific Computing Research, as part of the CHaRMNET Mathematical Multifaceted Integrated Capability Center (MMICC) program, under Award Number DE-SC0023164. IM release: LLNL-CONF-854965.
\end{ack}


\bibliography{ref}
\bibliographystyle{unsrt}

\appendix
\section*{Appendix: Additional Figures and Tables}
\setcounter{figure}{0}
\renewcommand{\thefigure}{A\arabic{figure}}
\setcounter{table}{0}
\renewcommand{\thetable}{A\arabic{table}}
\begin{table}[H]
    \centering
    \caption{Comparison of the number of trainable parameters and the training accuracy}
    \label{tab:train}
    \begin{tabular}{lccc}
        \toprule
        \textbf{Model} & \textbf{\# Parameters} & \textbf{Avg. $\frac{ \norm{\ubm_{\text{tr}} - \hat\ubm_{\text{tr}}} }{ \norm{\ubm_{\text{tr}}} }$} & \textbf{Max. $\frac{ \norm{\ubm_{\text{tr}} - \hat\ubm_{\text{tr}}} }{ \norm{\ubm_{\text{tr}}} }$} \\
        \midrule
        NODE         & 209128 & 1.798e-01 & 4.802e-01 \\
        NODE+PI      & 209128 & 1.798e-01 & 4.801e-01 \\
        PNODE        & 209641 & 1.094e-02 & 1.459e-02 \\
        PNODE+PI     & 209641 & 1.430e-02 & 1.931e-02 \\
        HyerPNODE    & 165509 & 8.336e-03 & 1.068e-02 \\
        HyerPNODE+PI & 165509 & 1.140e-02 & 1.421e-02 \\
        \bottomrule
    \end{tabular}
\end{table}

\begin{table}[H]
    \centering
    \caption{Comparison of the test accuracy and speedups}
    \label{tab:test}
    \begin{tabular}{lccccc}
        \toprule
        \textbf{Model} & \textbf{Avg. $\frac{ \norm{\ubm_{\text{ts}} - \hat\ubm_{\text{ts}}} }{ \norm{\ubm_{\text{ts}}} }$} & \textbf{Max. $\frac{ \norm{\ubm_{\text{ts}} - \hat\ubm_{\text{ts}}} }{ \norm{\ubm_{\text{ts}}} }$} & \textbf{Min. Speedup} & \textbf{Max. Speedup}\\
        \midrule
        NODE         & 2.653e-01 & 6.821e-01 & 1.133e+02 & 1.320e+03 \\
        NODE+PI      & 2.654e-01 & 6.819e-01 & 1.142e+02 & 1.364e+03 \\
        PNODE        & 3.212e-02 & 7.487e-02 & 1.164e+02 & 1.292e+03 \\
        PNODE+PI     & 4.861e-02 & 1.099e-01 & 1.151e+02 & 1.280e+03 \\
        HyerPNODE    & 4.537e-02 & 1.238e-01 & 1.173e+02 & 1.832e+03 \\
        HyerPNODE+PI & 4.618e-02 & 1.234e-01 & 1.167e+02 & 1.811e+03 \\
        \bottomrule
    \end{tabular}
\end{table}
\begin{figure}[H]
\centering
\begin{subfigure}{0.245\textwidth}
\centering
\includegraphics[width=\textwidth]{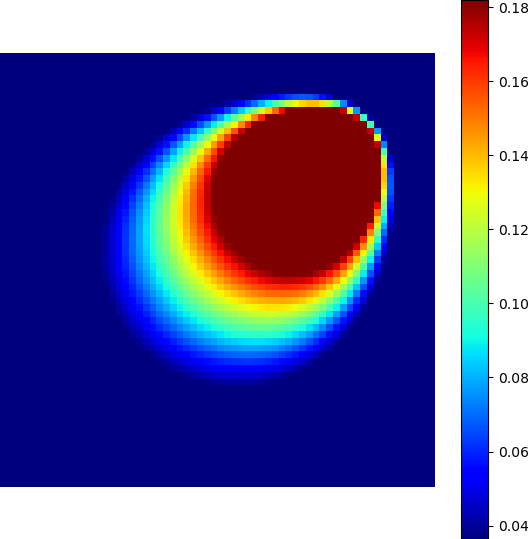}
\end{subfigure}
\hfill
\begin{subfigure}{0.245\textwidth}
\centering
\includegraphics[width=\textwidth]{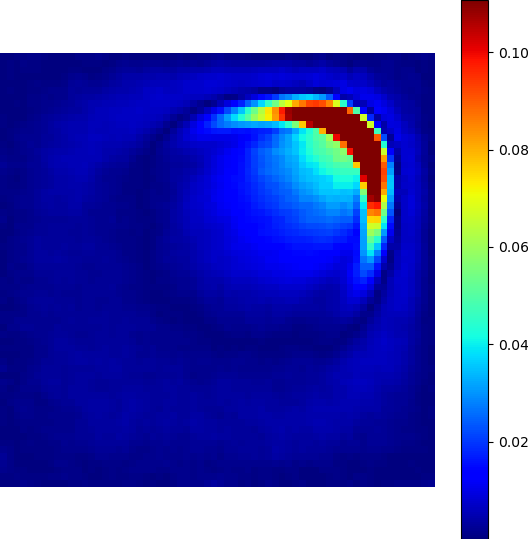}
\end{subfigure}
\hfill
\begin{subfigure}{0.245\textwidth}
\centering
\includegraphics[width=\textwidth]{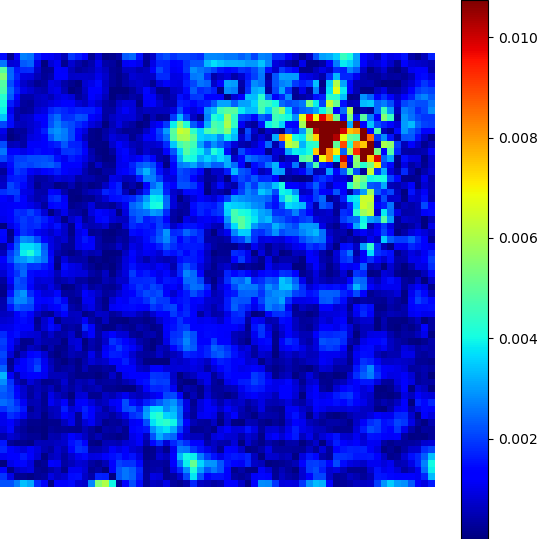}
\end{subfigure}
\hfill
\begin{subfigure}{0.245\textwidth}
\centering
\includegraphics[width=\textwidth]{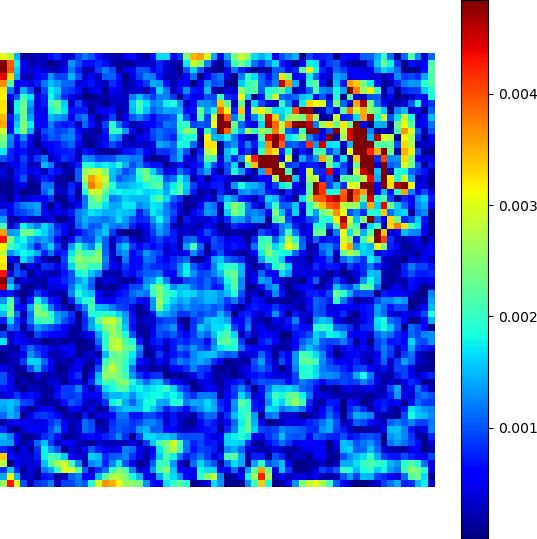}
\end{subfigure}
\caption{\small From left to right: FOM reference solution and the comparison of point-wise $\norm{\ubm_{\text{tr}} - \hat\ubm_{\text{tr}}}$ at $T=1$ and selected $Re = 500$ in the training set for NODE, PNODE, and HyperPNODE, respectively.
}
\label{fig:train}
\end{figure}
\begin{figure}[H]
\centering
\begin{subfigure}{0.245\textwidth}
\centering
\includegraphics[width=\textwidth]{_img/bg2d/node/bg2d_node_tr_state_0_reference.png}
\end{subfigure}
\hfill
\begin{subfigure}{0.245\textwidth}
\centering
\includegraphics[width=\textwidth]{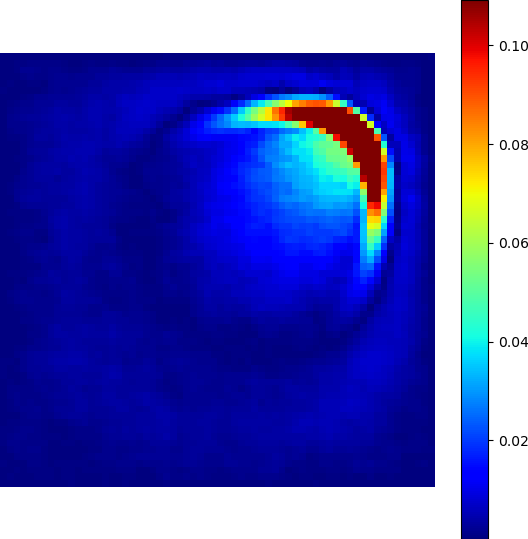}
\end{subfigure}
\hfill
\begin{subfigure}{0.245\textwidth}
\centering
\includegraphics[width=\textwidth]{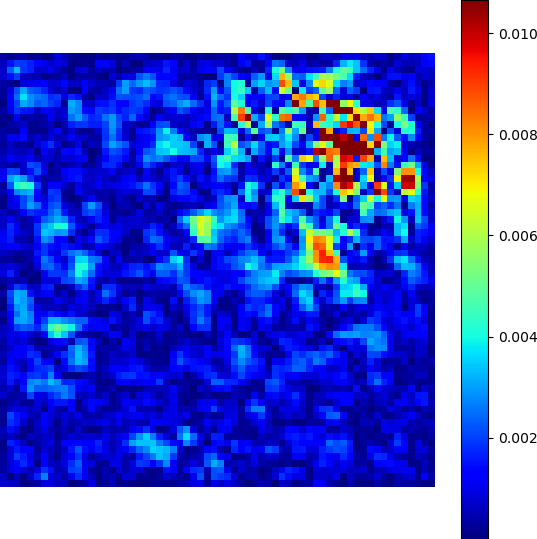}
\end{subfigure}
\hfill
\begin{subfigure}{0.245\textwidth}
\centering
\includegraphics[width=\textwidth]{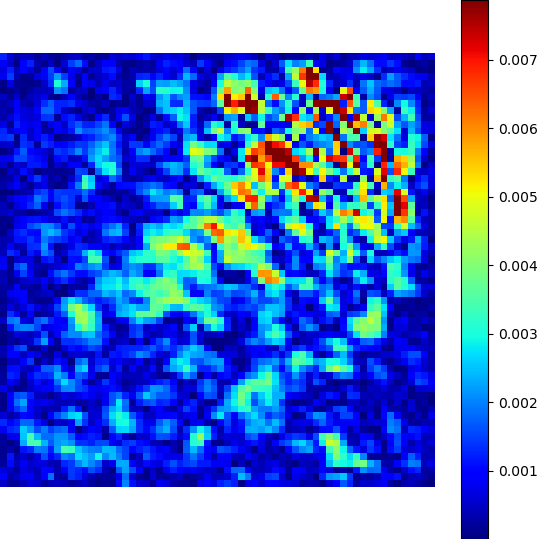}
\end{subfigure}
\caption{\small From left to right: FOM reference solution and the comparison of point-wise $\norm{\ubm_{\text{tr}} - \hat\ubm_{\text{tr}}}$ at $T=1$ and selected $Re = 500$ in the training set for NODE+PI, PNODE+PI, and HyperPNODE+PI, respectively.
}
\label{fig:train_pinn}
\end{figure}
\begin{figure}[H]
\centering
\begin{minipage}[t]{0.245\textwidth}
\centering
\includegraphics[width=1\textwidth]{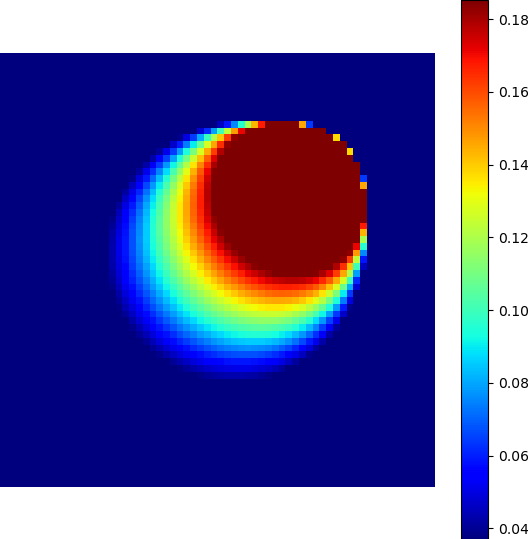}
\end{minipage}
\hfill
\begin{minipage}[t]{0.245\textwidth}
\centering
\includegraphics[width=\textwidth]{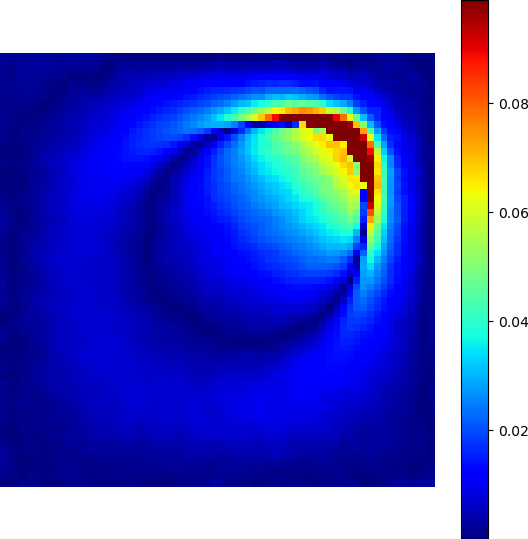}
\end{minipage}
\hfill
\begin{minipage}[t]{0.245\textwidth}
\centering
\includegraphics[width=\textwidth]{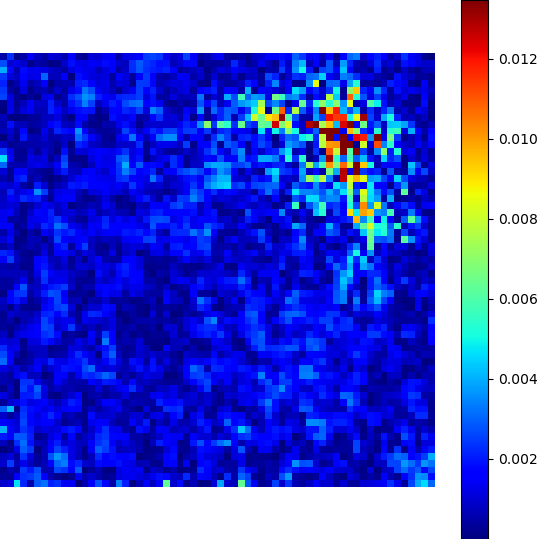}
\end{minipage}
\hfill
\begin{minipage}[t]{0.245\textwidth}
\centering
\includegraphics[width=\textwidth]{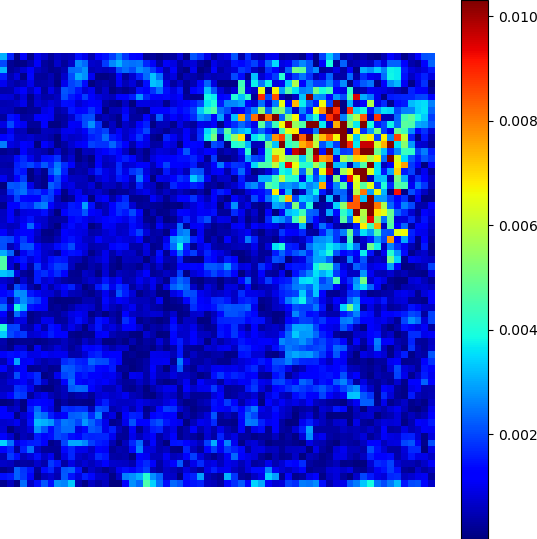}
\end{minipage}
\caption{\small From left to right: FOM reference solution and the comparison of point-wise $\norm{\ubm_{\text{ts}} - \hat\ubm_{\text{ts}}}$ at $T=1$ and selected $Re = 60000$ in the training set for NODE, PNODE, and HyperPNODE, respectively.
}
\label{fig:test}
\end{figure}
\begin{figure}[H]
\centering
\begin{minipage}[t]{0.245\textwidth}
\centering
\includegraphics[width=1\textwidth]{_img/bg2d/node/bg2d_node_ts_state_0_reference.png}
\end{minipage}
\hfill
\begin{minipage}[t]{0.245\textwidth}
\centering
\includegraphics[width=\textwidth]{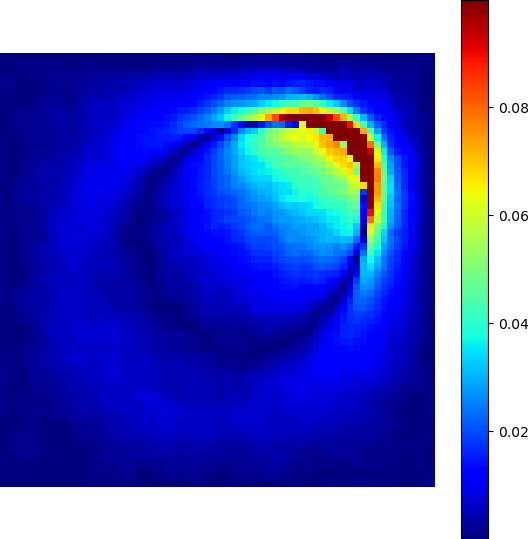}
\end{minipage}
\hfill
\begin{minipage}[t]{0.245\textwidth}
\centering
\includegraphics[width=\textwidth]{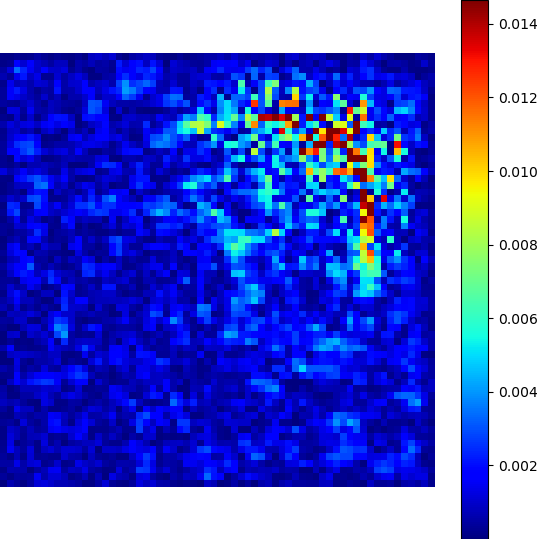}
\end{minipage}
\hfill
\begin{minipage}[t]{0.245\textwidth}
\centering
\includegraphics[width=\textwidth]{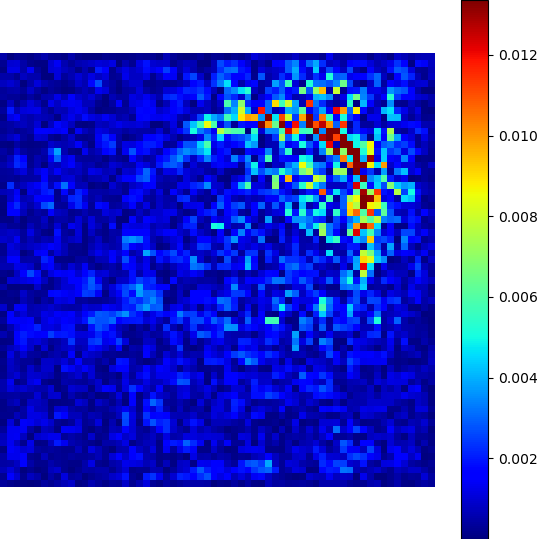}
\end{minipage}
\caption{\small From left to right: FOM reference solution and the comparison of point-wise $\norm{\ubm_{\text{ts}} - \hat\ubm_{\text{ts}}}$ at $T=1$ and selected $Re = 60000$ in the training set for NODE+PI, PNODE+PI, and HyperPNODE+PI, respectively.
}
\label{fig:test_pinn}
\end{figure}

\clearpage


\end{document}